\documentclass[11pt]{article}

\usepackage{amsmath}
\usepackage{amssymb}
\usepackage{eucal}

\begin{document}

\baselineskip 16pt

\title{On the   lattice
 of the $\sigma$-permutable subgroups of a finite group}

\author            
{ Alexander  N. Skiba \\
{\small Department of Mathematics and Technologies of Programming,}\\
{\small Francisk Skorina Gomel State University,}\\
{\small Gomel 246019, Belarus}\\
{\small E-mail: alexander.skiba49@gmail.com}}

\date{}
\maketitle

\begin{abstract}    Let   $\sigma =\{\sigma_{i} | i\in I\}$ be some
 partition of the set of all primes $\Bbb{P}$, $G$ a  finite group and $\sigma (G)
=\{\sigma_{i} |\sigma_{i}\cap \pi (G)\ne  \emptyset \}$. 

A set
 ${\cal H}$ of subgroups of $G$  is said to be a  \emph{complete
Hall $\sigma $-set} of $G$  if every member $\ne 1$ of  ${\cal H}$
 is a Hall $\sigma _{i}$-subgroup of $G$ for some $\sigma _{i}\in \sigma $ and
   ${\cal H}$ contains exactly one Hall  $\sigma _{i}$-subgroup of $G$ for every
 $\sigma _{i}\in  \sigma (G)$.
A subgroup $A$ of $G$ is said to be  \emph{${\sigma}$-permutable}   in $G$ if $G$ possesses a complete
Hall $\sigma $-set  and   $A$ permutes with each   Hall $\sigma _{i}$-subgroup $H$ of
 $G$, that is,  $AH=HA$ for all $i \in I$.

We characterize finite groups with  distributive lattice of 
the ${\sigma}$-permutable subgroups.

\end{abstract}                                      


\let\thefootnote\thefootnoteorig

\section{Introduction}

\

Throughout this paper,  $G$ always denotes
a finite group. Moreover, we  
use ${\cal L}(G)$ to denote the lattice of all subgroups of $G$.
 $\mathbb{P}$ is the set of all  primes,
  $\pi \subseteq  \Bbb{P}$ and  $\pi' =  \Bbb{P} \setminus \pi$.
As usual,  $\pi (G)$ is  the set of all
  primes dividing the order $|G|$ of $G$.   
The subgroups $A$ and $B$ of $G$   are said to be   
\emph{permutable}  if $AB=BA$.  In this case they also say that $A$ 
\emph{permutes } with $B$. If $A$      permutes with all Sylow subgroups 
of $G$, then $A$ is called  \emph{$S$-permutable} in $G$ \cite{prod}.
Recall also that an element  $a $     of the  lattice $\cal L$  is called 
\emph{meet-distributive} \cite[p. 136]{R.}
 if $a  \wedge (b   \vee c) =(a \wedge b) \vee  (a    \wedge c)$
for all  $b, c\in \cal L$.

In what follows, $\sigma =\{\sigma_{i} | i\in I\}$ is 
 some
partition of $\Bbb{P}$, that is,
$\Bbb{P}=\bigcup_{i\in I} \sigma_{i}$ and $\sigma_{i}\cap
\sigma_{j}= \emptyset $ for all $i\ne j$.

A  set  ${\cal H}$ of subgroups of $G$ is a
 \emph{complete Hall $\sigma $-set} of $G$ \cite{3}  if
 every member $\ne 1$ of  ${\cal H}$ is a Hall $\sigma _{i}$-subgroup of $G$
 for some $\sigma _{i} \in \sigma$ and ${\cal H}$ contains exactly one Hall
 $\sigma _{i}$-subgroup of $G$ for every $i$ such that   $\sigma _{i}\cap  \pi (G)\ne \emptyset$.  
$G$ is said   to be   \emph{ $\sigma $-full} if it possesses a complete 
Hall $\sigma$-set.

{\bf Definition 1.1.}  We say  that  a
  subgroup $A$ of $G$ is said to be  \emph{$\sigma$-quasinormal} or
    \emph{$\sigma$-permutable}  in $G$  \cite{1} if $G$ is $\sigma$-full and $A$
 permutes with each   Hall $\sigma _{i}$-subgroup $H$ of
 $G$  for all $i \in I$.

{\bf Remark 1.2. }   (i) If $G$ possesses
 a complete Hall $\sigma$-set  ${\cal H}$ such that $AH^{x}=H^{x}A$ for all
  $H\in {\cal H}$ and  all $x\in G$, then  $A$ is $\sigma$-permutable   in $G$
 (see Proposition 3.1 below).

(ii)  $G$ is called \emph{$\sigma$-decomposable} (Shemetkov \cite{Shem})  or 
\emph{$\sigma$-nilpotent} (Guo and Skiba  \cite{33}) if $G=H_{1}\times \cdots \times H_{t}$,
 where 
$\{H_{1}, \ldots ,   H_{t}\}$ is a complete Hall $\sigma$-set of $G$.  It 
is not difficult to show that $G$ is $\sigma$-nilpotent if and only if 
every subgroup of $G$ is $\sigma$-permutable in $G$.

(iii) In the classical case
 when $\sigma =\sigma ^{0}=\{\{2\}, \{3\}, \ldots 
\}$: $G$  is  $\sigma ^{0}$-nilpotent if
 and only if $G$ is  nilpotent; a subgroup $A$ of $G$ is $\sigma^{0}$-permutable
 in $G$ if and only if it is $S$-permutable in $G$.

(iv) In the other classical case when  $\sigma =\sigma ^{\pi}=\{\pi, \pi'\}$:
 $G$  is  $\sigma ^{\pi}$-nilpotent 
 if and only if it is $\pi$-decomposable, that is, 
  $G=O_{\pi}(G)\times O_{\pi'}(G)$; a subgroup $A$ of a  $\pi$-separable group
 $G$ is $\sigma^{\pi}$-permutable in $G$ if and only if $A$ permutes  with 
all Hall $\pi$-subgroups  and with all Hall $\pi'$-subgroups of $G$.

(v) If $G$ is $\sigma$-full,  the set  ${\cal L}_{\sigma per}(G)$  of all $\sigma$-permutable 
subgroups of  $G$   
is partially ordered with respect to   set inclusion. Moreover, ${\cal 
L}_{\sigma per}(G)$   is a lattice  
 since $1\in {\cal L}_{\sigma per}(G)$ and,  by Lemma 2.1 below,   for any 
 $A_{1}, \ldots , A_{n}\in {\cal L}_{\sigma per}(G)$ the subgroup $\langle
 A_{1}, \ldots , A_{n}\rangle $ 
  is the least upper bound for $\{A_{1}, \ldots , A_{n}\}$  in  ${\cal 
L}_{\sigma per}(G)$.

The conditions under which the lattice  ${\cal L}_{sn}(G)$ of all subnormal 
subgroups of $G$   is modular or 
distributive are known (see \cite[Theorems 9.2.3, 9.2.4]{R.}).
It is well-known also  that the lattice ${\cal L}_{n}(G)$ of all  normal 
subgroups of $G$ is modular and this lattice is  distributive if and only if
in every factor group $G/R$, any two $G/R$-isomorphic normal subgroups coincide
   (see \cite{Pazd} and   \cite[Theorem 9.1.6]{R.}).  Kegel proved \cite{Keg} that the  
set ${\cal L}_{S}(G)$  of all $S$-permutable subgroups of $G$  forms  a sublattice 
of the lattice ${\cal L}_{sn}(G)$.  Since  $ {\cal L}_{n}(G)  \subseteq 
 {\cal L}_{S}(G) \subseteq {\cal L}_{sn}(G) $, where  both  
inclusions in general  are strict, it seems natural to ask:  {\sl  Under what conditions 
 the lattice ${\cal L}_{S}(G) $
 is  modular or distributive?}  
 Moreover, in view of Remark 1.2(v), it makes sense to consider  the 
following general

{\bf Question 1.3} (See Questions 6.10 and 6.11 in \cite{3}). {\sl  Under what conditions 
 the lattice ${\cal L}_{\sigma per}(G) $
 is  modular or distributive?}

 Note that if  $K\trianglelefteq H$  and $K, H \in {\cal L}_{\sigma _{i}}(G)$, where  
${\cal L}_{\sigma _{i}}(G)$ is the set of all $\sigma$-permutable $\sigma _{i}$-subgroups
 of $G$,
 then $O^{\sigma _{i}}(G)$ normalizes both subgroups $K$ and $H$ (see Lemma
 2.4(1) below) and hence we can consider  
$O^{\sigma _{i}}(G)$ as a  group of operators for $H/K$ (assuming, as usual,  that
 $(hK)^{a}=h^{a}K$ for all $hK\in H/K$  and 
 $a\in O^{\sigma _{i}}(G)$).

We do not know under what conditions on $G$  the lattice ${\cal L}_{\sigma}(G)$  is modular.
Nevertheless, we give a full answer to the second part of Question 1.3.

{\bf Theorem A.} {\sl   Suppose that  $G$ is $\sigma$-full. Let
  $D=G^{{\frak{N}}_{\sigma}} $ and  ${\cal L}={\cal L}_{\sigma}(G)$. 
Then the lattice ${\cal L}$ is distributive if and only if
 the following conditions hold:}

(i) {\sl Every two members 
 of ${\cal L}$  are  permutable.}

(ii) {\sl The lattice ${\cal L}_{n}(G)$ of all  normal subgroups of $G$ is distributive. }

(iii) {\sl $G/D$ is cyclic and $D$ is a  meet-distributive element of   ${\cal L}$.}

(iv) {\sl  In every factor group $\bar{G}=G/R$,   any 
 two $O^{\sigma _{i}}(\bar{G})$-isomorphic sections
 $\bar{H}/\bar{K}$ and $\bar{L}/\bar{K}$, where   $\bar{K}, \bar{H}, 
\bar{L}\in {\cal L}_{\sigma _{i}}(\bar{G})$ for some $i$, coincide.}

In this theorem  $G^{{\frak{N}}_{\sigma}}$  denotes the \emph{$\sigma$-nilpotent residual}
 of $G$,  that is, the intersection of all normal subgroups $N$ of $G$ 
with $\sigma$-nilpotent quotient $G/N$. 

 Theorem A remains to be new for each special  partition of $\Bbb{P}$.   
In particular, in the case when $\sigma =\sigma ^{0}$ we get from Theorem A  the 
following 

{\bf Corollary  1.4.}
 {\sl   Let  $D=G^{\frak{N}}$ be the
 nilpotent residual of $G$
 and  ${\cal L}={\cal 
 a L}_{S}(G)$.   Then the lattice ${\cal L}$ is distributive   if and only if
 the following conditions hold:}

(1) {\sl Conditions (i), (ii) and (iii) in Theorem A fold  for $G$. }

(2) {\sl  In every factor group $\bar{G}=G/R$, any two $O^{p}(\bar{G})$-isomorphic sections
 $\bar{H}/\bar{K}$ and $\bar{L}/\bar{K}$, where   $\bar{K}, \bar{H}, 
\bar{L}\in {\cal L}_{pS}(\bar{G})$ and $p$ is a prime, coincide.}

In this corollary ${\cal L}_{pS}(\bar{G})$ denotes  the set of all
 $S$-permutable $p$-subgroups of  $\bar{G}$.

In the case when  $\sigma =\sigma ^{\pi}$ (see Remark 1.2(iv)) we get from 
Theorem A the following fact.

{\bf Corollary  1.5.} {\sl   Let  $D$ be the
 $\pi$-decomposable  residual of $G$, that is, the smallest normal subgroup 
of $G$ with  $\pi$-decomposable quotient $G/D$. Suppose that $G$ is $\pi$-separable and let 
${\cal L}={\cal L}_{\sigma ^{\pi}}(G)$.  Then the lattice
 ${\cal L}$ is distributive   if and only if the following conditions hold:}

(1) {\sl Conditions (i), (ii) and (iii) in Theorem A hold  for $G$. }

(2) {\sl  In every factor group $\bar{G}=G/R$, any two $O^{\pi}(\bar{G})$-isomorphic sections
 $\bar{H}/\bar{K}$ and $\bar{L}/\bar{K}$, where 
 $\bar{K}, \bar{H}$ and  
$ \bar{L}$ are  $\sigma ^{\pi}$-permutable $\pi$-subgroups of $\bar{G}$, coincide.}

(3)  {\sl  In every factor group $\bar{G}=G/R$, any two $O^{\pi'}(\bar{G})$-isomorphic sections
 $\bar{H}/\bar{K}$ and $\bar{L}/\bar{K}$, where 
 $\bar{K}, \bar{H}$ and  
$\bar{L}$ are  $\sigma ^{\pi}$-permutable $\pi'$-subgroups of $\bar{G}$, coincide.}

The proof of Theorem A consists of many steps and the next theorems are 
two  of them.

{\bf Theorem B.} {\sl  Suppose that  $G$ is $\sigma$-full. Then ${\cal L}_{\sigma per}(G)$ is 
 a sublattice of the lattice ${\cal L}(G)$.}

{\bf Corollary 1.6} (Kegel \cite{Keg}).  {\sl  The set
$ {\cal L}_{S}(G)$ of all
$S$-permutable subgroups of $G$ forms a sublattice of the lattice ${\cal L}(G)$.}

There are at least three different proofs of Corollary 1.6 (see, for 
example,  \cite{Keg, 99,  1}). One more,  the shortest one,    
gives the proof of Theorem B.

{\bf Theorem  C.}  {\sl A $\sigma$-nilpotent subgroup $A$ of 
 $G$ is $\sigma$-permutable in $G$ if and only if each
 characteristic subgroup of $A$  is $\sigma$-permutable 
in $G$. }

{\bf Corollary 1.7} (See \cite{99} or \cite[Theorem 1.2.17]{prod}).  {\sl Let $A$ be
 a nilpotent subgroup of $G$.
  Then the following   statements are equivalent:  }

(i) {\sl $A$ is $S$-permutable in $G$.}

(ii) {\sl Each Sylow  subgroup of $A$  is $S$-permutable 
in $G$.  }

(iii) {\sl Each characteristic subgroup of $A$  is $S$-permutable 
in $G$. }

\section{Proof of Theorems B  and C}

{\bf Lemma  2.1} (See  \cite[A, Lemma 1.6]{DH}). {\sl Let  $A$,  $B$ and $H$ be subgroups of $G$.
If $AH=HA$ and $BH=HB$, then $\langle A, B\rangle H=H\langle A, B\rangle$.    }

A  subgroup $A$ of $G$ is called \emph{${\sigma}$-subnormal} in $G$
 \cite{1}  if there is a subgroup chain  $A=A_{0} \leq A_{1} \leq \cdots \leq
A_{t}=G$  such that  either $A_{i-1}\trianglelefteq A_{i}$ or $A_{i}/(A_{i-1})_{A_{i}}$
 is  ${\sigma}$-primary for all $i=1, \ldots , t$.

The importance of this  concept  is related to the 
following result.
 
{\bf Lemma  2.2} (See  \cite[Theorem B]{1}). {\sl Let  $A$  be a subgroup of $G$. If $G$ 
 possesses  
 a complete Hall $\sigma$-set  ${\cal H}$ such that $AH^{x}=H^{x}A$ for all
  $H\in {\cal H}$ and  all $x\in G$, then  $A$ is  $\sigma$-subnormal in $G$.    }

{\bf Lemma  2.3} (See Lemma 2.6 in \cite{1}). {\sl Let  $A$,  $K$ and $N$ be subgroups of 
  $G$.
 Suppose that   $A$
is $\sigma $-subnormal  in $G$ and $N$ is normal in $G$.
    }

(1) {\sl $A\cap K$    is  $\sigma $-subnormal in
$K$}.

(2) {\sl If  $|G:A|$ is a $\sigma _{i}$-number,  then 
 $O^{\sigma _{i}}(A)= O^{\sigma _{i}}(G)$.}

(3) {\sl $AN/N$ is
$\sigma$-subnormal in $G/N$. }

(4) {\sl If    $A $
  is a $\sigma _{i}$-group,  then  $A\leq O_{\sigma _{i}}(G)$. }

(5) {\sl If $H\ne 1 $ is a Hall $\sigma _{i}$-subgroup of $G$, then $A\cap H\ne 1$ is
 a Hall $\sigma _{i}$-subgroup of $A$. }

(6) {\sl  If $K$ is  $\sigma $-subnormal  in $G$ and the subgroups  $A$ and $K$ are
 $\sigma$-nilpotent, then $\langle A, K\rangle$  is $\sigma$-nilpotent.}

{\bf Proof of Theorem B.}  
 In fact, in view of Lemmas 2.1 and 2.2, it is enough  to show
 that if $A$ and  $B$ are $\sigma$-subnormal subgroups  of $G$ such that for a Hall
 $\sigma _{i}$-subgroup $H$ of $G$ we have $AH=HA$ and $BH=HB$, then $(A\cap B)H=H(A\cap B)$. 
 Assume that this is false and let $G$ be a 
counterexample of minimal order. Then $G$ is not a $\sigma _{i}$-group, since otherwise 
we have $H=G$ and so $G=(A\cap B)H=H(A\cap B)$.

 Let  $E=AH\cap BH$. Then $A\cap E$ and $B\cap E$ are $\sigma$-subnormal subgroups 
 in 
$E$  by  Lemma 2.3(1). Moreover, $AH\cap E=H(A\cap E)=(A\cap E)H$. 
Similarly,  $(B\cap E)H=H(B\cap E)$. Hence the hypothesis holds for
 $(A\cap E, B\cap E, H, E)$.  
 Assume that $E < G$. Then the choice of $G$ implies that
 $A\cap B=(A\cap E)\cap (B\cap E)$ is  permutable with $H$.
 Hence   $E=G$, 
so $G=AH=BH$. Thus $|G:A|$ and $|G:B|$ are  $\sigma _{i}$-numbers. 
  Hence by 
Lemma 2.3(2)  we have $O^{\sigma _{i}}(A)=O^{\sigma _{i}}(G)=O^{\sigma _{i}}(B).$
Therefore, since $G$ is not a $\sigma _{i}$-group, it follows that  $V=A_{G}\cap 
B_{G}\ne 1$.  Moreover, $A/V$ and $B/V$ are  $\sigma$-subnormal subgroups 
of $G/V$ by Lemma 2.3(3). Also we have $$(A/V)(HV/V)=AH/V=HA/V=(HV/V)(A/V)$$ 
and  $(B/V)(HV/V)=(HV/V)(B/V)$, where $HV/V$ is 
a Hall $\sigma _{i}$-subgroup of $G/V$.
Hence the choice of $G$ implies that $$(A\cap B/V)(HV/V)=   ((A/V)\cap
 (B/V))(HV/V)=$$$$=(HV/V)((A/V)\cap (B/V))
=(HV/V)(A\cap B/V).$$  But then  $$(A\cap B)H=(A\cap B)HV=HV(A\cap B)= 
H(A\cap B).$$
 This contradiction completes the proof of the result.

{\bf Lemma 2.4}  (See Lemmas 2.8,   3.1 and Theorem B in \cite{1}).  {\sl 
Let  $A$ and $B$ be subgroups of $G$. Suppose that $G$ 
 possesses  
 a complete Hall $\sigma$-set  ${\cal H}$ such that $AH^{x}=H^{x}A$ for all
  $H\in {\cal H}$ and  all $x\in G$.  Then:}

(1) {\sl If $A$ is a $\sigma _{i}$-group, then $O^{\sigma _{i}}(G)  \leq N_{G}(A)$}.

(2)  {\sl  $A/A_{G}$ is $\sigma$-nilpotent}.

(3)  {\sl If $B$ is a $\sigma _{i}$-group and  $O^{\sigma _{i}}(G)  \leq 
N_{G}(B)$, then $G$ 
 possesses  
 a complete Hall $\sigma$-set  ${\cal L}$ such that $BL^{x}=L^{x}B$ for all
  $L\in {\cal L}$ and  all $x\in G$}.

{\bf  Proposition 2.5.} {\sl  Let   $A$ be a  $\sigma$-nilpotent 
 $\sigma$-subnormal  subgroup of $G$ and  $V$  a characteristic subgroup 
of $A$. Let  $H$ be a Hall $\sigma _{i}$-subgroup of $G$. 
If $AH=HA$, then  $VH=HV$.}

{\bf Proof. }    Assume that this proposition is false and let $G$ be a counterexample 
with   $|G|+|V|+|A|$   minimal.

By hypothesis,  
$A=A_{1}\times \cdots \times A_{t}$,
 where $\{A_{1}, \ldots ,  A_{t}\}$ is a
 complete  Hall $\sigma$-set of $A$. Hence  $V=(A_{1}\cap V) \times \cdots \times (A_{t} 
\cap V)$,
 where $\{A_{1}\cap V, \ldots ,  A_{t}\cap V\}$ is a
 complete  Hall $\sigma$-set of $V$.  We can assume without loss of 
generality that $A_{k}$ is a   $\sigma _{k}$-subgroup of $A$  for all 
$k=1, \ldots, t$.                  

 It is clear that $A_{i}\cap V$ is 
 characteristic in $A$  for all $i=1, \ldots, t$.  Therefore,  if 
$A_{i}\cap V < V$, then    $(A_{i}\cap V)H=H(A_{i}\cap V)$ by the choice 
of $G$ and  so   for some $j$, $j=1$  say,  we have  $A_{1}\cap V = V$
 since otherwise we have $$VH=((A_{1}\cap V) \times \cdots \times (A_{t} 
\cap V))H=H((A_{1}\cap V) \times \cdots \times (A_{t} 
\cap V))=HV.$$ Thus $V\leq A_{1}$. It is clear that  $A_{1}$ is
 a $\sigma$-subnormal subgroup of $G$, so in the case when $i=1$
  we have $V\leq A_{1}\leq H$ by Lemma 2.3(5). But then
$VH=H=HV$, a contradiction. Thus $ i > 1.$ 

Now  we  show that $A_{1}H=HA_{1}$. 
 Indeed,  it is clear that $A=A_{1}\times V\times A_{i}$ and   $A_{i}$ is 
$\sigma$-subnormal in $G$. Thus  $A_{i}\leq  H$ by Lemma 2.3(5). Therefore 
 $$AH=HA=(A_{1}\times V\times A_{i})H=(A_{1}\times V)H=(A_{1}\times V)H,$$
 where $A_{1}\times V$ is
 a $\sigma$-subnormal 
$\sigma _{i}'$-subgroup of $G$. Then $A_{1}\times V$ is
  $\sigma$-subnormal  in $(A_{1}\times V)H$ by Lemma 2.3(1).
   Hence $H\leq N_{G}(A_{1}\times V)$ by Lemma 2.3(2). Since 
$A_{1}$ is a characteristic subgroup of $A_{1}\times V$, we have $H\leq 
N_{G}(A_{1})$  and so   $A_{1}H=HA_{1}$.  Therefore $H\leq N_{G}(A_{1})$ 
by Lemma 2.3(2). But  $V$ is a characteristic subgroup 
of $A_{1}$ since  $V$ is  characteristic  in $A$ by hypothesis and
 $A=A_{1}\times \cdots \times A_{t}$.  Therefore  $H\leq N_{G}(V)$ and so 
$VH=HV$, a contradiction.  
 The proposition  is proved.

{\bf Corollary 2.6.}  {\sl Let $A$ be a $\sigma$-nilpotent subgroup of a
  $\sigma$-full group $G$.
  Then the following   statements are equivalent:  }

(i) {\sl $A$ is $\sigma$-permutable in $G$.}

(ii) {\sl Each Hall $\sigma _{i}$-subgroup of $A$  is $\sigma$-permutable 
in $G$ for all $i$.    }

(iii) {\sl Each characteristic subgroup of $A$  is $\sigma$-permutable 
in $G$. }

{\bf Proof. }  By hypothesis, $A=A_{1}\times \cdots \times A_{t}$, where
 $\{A_{1}, \ldots ,  A_{t}\}$ is a
 complete  Hall $\sigma$-set of $A$. Then $A_{i}$ is characteristic in $A$ 
for all $i=1, \ldots ,t$.   Therefore 
(ii), (iii) $\Rightarrow $ (i).

 (i) $\Rightarrow$ (ii), (iii)  This  follows from 
Proposition 2.5.

 The corollary   is proved. 

{\bf Proof of Theorem C. }  This directly follows from Corollary 2.6.

 \section{Proof of Theorem A}

{\bf   Proposition 3.1.} {\sl  Let $A$ be a subgroup of $G$. If $G$ 
 possesses  
 a complete Hall $\sigma$-set  ${\cal H}$ such that $AL^{x}=L^{x}A$ for all
  $L\in {\cal H}$ and  all $x\in G$, then  $A$ 
  is $\sigma$-permutable in $G$. }

 {\bf Proof. }  Assume that this proposition is false and let $G$ be a counterexample 
with   $|G|+|A|$   minimal.  Then for some  $i$ and some Hall $\sigma _{i}$-subgroup 
$H$ of $G$ we have $AH\ne HA$.  Let  ${\cal H}=\{H_{1},\ldots , H_{t}\}$. 
We  can  assume without loss of generality that $H_{k}$ is a $\sigma 
_{k}$-group for all $k=1,  \ldots , t$.  Let $V=H_{i}$.

First we show that $A_{G}=1$. Indeed, assume that $R=A_{G}\ne 1$.  
Then  ${\cal H}_{0}=\{H_{1}R/R,\ldots , H_{t}R/R\}$ is a complete Hall $\sigma$-set
 of  $G/R$ such that   $$AL^{x}/R=(A/R)(LR/R)^{xR}=(LR/R)^{xR}(A/R)=L^{x}A/R$$ for all
  $LR/R\in {\cal H}_{0}$ and  all $xR\in G/R$. On the other hand, $HR/R$   is Hall
 $\sigma _{i}$-subgroup of $G/R$.  
Hence the choice of $G$ implies that $$AH/R=(A/R)(HR/R)=(HR/R)(A/R)=HA/R$$ and
  so $AH= HA$, a contradiction.    
Therefore $A_{G}=1$, hence $A=A_{1}\times \cdots \times A_{t}$,
 where $\{A_{1}, \ldots ,  A_{t}\}$ is a
 complete  Hall $\sigma$-set of $A$ by Lemma 2.4(2).  Moreover, Lemma 2.2 
implies that $A$ is $\sigma$-subnormal in $G$.

 First assume that $A=A_{1}$ is a 
$\sigma _{j}$-group.  If $j=i$, then  $A\cap H=A$ by Lemma 2.3(5) and  so 
$AH=H=HA$. Hence $j\ne i$.   
 By hypothesis,   $AV^{x}=V^{x}A$ for each   $x\in G$.  Then   $V^{x}\leq N_{G}(A)$
 for all  $x\in G$ by Lemma 
2.3(1)(2). Hence $V^{G}\leq N_{G}(A)$. But then $H\leq   V^{G}\leq 
N_{G}(A)$, which implies that $AH=HA$.  This contradiction shows that   $ 
t > 1$.

The subgroups $A_{1}, \ldots ,  A_{t}$ are characteristic in $A$, so   $A_{i}L^{x}=L^{x}A_{i}$
 for all
  $L\in {\cal H}$ and  all $x\in G$  by Proposition 2.5.
Therefore the minimality of $|G|+|A|$   implies that 
$A_{i}H=HA_{i}$     for all $i=1, \ldots , t$, so  $AH=HA$. 
  This contradiction 
completes the proof of the result. 

{\bf Lemma 3.2.}  {\sl Let  $R\leq V$ and  $H$  be subgroups of a $\sigma$-full group  $G$.
 Suppose that  $H$ is
 $\sigma$-permutable in $G$ and $R$ is normal in $G$. Then:}

(1) {\sl  If $V/R$ is a 
$\sigma$-permutable subgroup of $G/R$, then $V$ is a $\sigma$-permutable 
subgroup of $G$.  }

(2)    {\sl The subgroup  $HR/R$ is $\sigma$-permutable in $G/R$.}

{\bf Proof. }   
(1) Let $i\in I$ and $H$ be  a Hall $\sigma _{i}$-subgroup of $G$.  Then 
$HR/R$ is  a  Hall $\sigma _{i}$-subgroup of $G/R$, so  $$VH/R=(V/R)(HR/R)=(HR/R)(V/R)=HV/R$$ 
by hypothesis and hence  $VH=HV$.

(2) By hypothesis,  $G$ possesses
 a complete Hall $\sigma$-set  ${\cal H}=\{H_{1},\ldots , H_{t}\}$ and $HL^{x}=L^{x}H$ for all
  $L\in {\cal H}$ and  all $x\in G$.   
  Then ${\cal H}_{0}=\{H_{1}R/R,\ldots , H_{t}R/R\}$ is a complete Hall $\sigma$-set
 in $G/R$ and     $$(HR/R)(LR/R)^{xR}=HL^{x}R/R=L^{x}HR/R=(LR/R)^{xR}(HR/R)$$ for all
  $LR/R\in {\cal H}_{0}$ and $xR\in G/R$. Therefore   $HR/R$ is 
$\sigma$-permutable   in $G/R$ by Proposition 3.1.

The lemma is proved.

{\bf  Lemma 3.3} (See Lemma 5.2 in \cite{kimber}).  {\sl Let $\cal L$ be 
 a modular sublattice of the lattice ${\cal L}(G)$, and $U, V , N \in {\cal L}$
 with $N \trianglelefteq  \langle U, V \rangle $.  If $U$  permutes both 
with 
  $V\cap UN$ and $VN$, then $U$  permutes with $V$.   }

{\bf Proposition 3.4.} {\sl  Let $G$ be   $\sigma$-full and    ${\cal L}={\cal L}_{\sigma _{i}}(G)$. 
Then:  (i)     ${\cal L}$  is a sublattice of ${\cal L}_{\sigma per}(G)$, and 
(ii)  If ${\cal L}$  
 is distributive, then  $AB=BA$ for all $A, B\in {\cal L}$.
  }

{\bf Proof.} (i)  Let $A, B\in {\cal L}$.   The subgroups $A$ and $B$ are $\sigma$-subnormal 
in $G$
 by Lemma 2.2, so $A, B\leq  O_{\sigma _{i}}(G)$ by Lemma 2.3(4).
 Thus $\langle A, B \rangle $
 is a $\sigma _{i}$-subgroup of $G$ and this subgroup is $\sigma$-permutable in
 $G$ by Lemma 2.1.
 Finally, $A\cap B$ is also a  $\sigma _{i}$-subgroup of $G$ and this subgroup is
 $\sigma$-permutable in $G$
 by Theorem B. Thus we have (i).

(ii) Suppose that this assertion  is false and  let $G$ be 
 a counterexample with $|G|+|A|+|B|$   minimal.  Thus $AB\ne BA$ but  
 $A_{1}B_{1}= B_{1}A_{1}$ for all $A_{1}, B_{1} \in{\cal L}$ such that $A_{1} 
\leq A$, $B_{1}\leq B$ and either  $A_{1}\ne A$ or $B_{1}\ne  B$. 
Let    $V=\langle A, B \rangle O^{\sigma _{i}}(G)$ and $R=\langle A, B 
\rangle\cap O^{\sigma _{i}}(G)$.  Then $V$ is $\sigma$-subnormal in $G$.

(1) {\sl    ${\cal L}_{\sigma _{i}}(V)$   is a sublattice 
 of ${\cal L}$.    }
  
Indeed, let $H\in {\cal L}_{\sigma _{i}}(V)$. Then $H$ is $\sigma$-subnormal in $V$ by 
Lemma 2.2 and  so, because of Lemma 2.3(4),   
  $H \leq O_{\sigma _{i}}(V)\leq O_{\sigma _{i} }(G)$. Therefore $H$ 
permutes with each Hall  $\sigma _{i}$-subgroup of $G$. On the other hand, 
each  Hall   $\sigma _{j}$-subgroup $W$ of $G$, where $j\ne i$,  is contained in $V$ by
 Lemma 2.3(5) since $|W|$ divides $|V|$, so   $HW=WH$.  Hence $H\in {\cal L}_{\sigma}(G)$, which implies  (1). 

(2) {\sl $V=G, $   so $\langle A, B \rangle \trianglelefteq G$. }

 Claim (1) implies that 
 the hypothesis holds for  ${\cal L}_{\sigma _{i}}(V)$ and so in the case when $V\ne G$ the
 choice of $G$ implies  
  that $AB=BA$. Thus  
 $G= \langle A, B \rangle O^{\sigma _{i}}(G)$.   
 Therefore, since  
$O^{\sigma _{i}}(G) \leq N_{G}(\langle A, B \rangle )$ by Lemma 2.4(1), 
 $\langle A, B \rangle$ is normal in $G$.

(3)  $R=1$.

  Assume that  $R=\langle A, B  \rangle\cap O^{\sigma _{i}}(G)   \ne 1$. 
First we show that    $BRA=\langle A, B \rangle R$. Indeed, let $H/R$ be a
  $\sigma _{i}$-subgroup of $G/R$.
  Then  $H$ is a  $\sigma _{i}$-group   
  since $\langle A, B \rangle \leq O_{\sigma _{i}}(G)$. Moreover, 
 Lemma 3.2(1)(2) implies that   $H/R$ 
is $\sigma$-permutable in $G/R$  if and only if 
$H$ is $\sigma$-permutable in $G$. Therefore the
 lattice ${\cal L}_{\sigma _{i}}(G/R)$ is  isomorphic
 to the interval $[G/R]$ in the distributive 
lattice ${\cal L}$. Therefore, by the minimality of $G$, 
$(AR/R)(BR/R)=(BR/R)(AR/R)$ by Lemma 3.2(2) and so $BRA=\langle A, B \rangle R$. 

Now we show that $BRA=BR$. Assume that this is false. Then  $A\cap BR < A$.
 But Theorem B implies 
that $A\cap BR $ is $\sigma$-permutable in $G$, so the minimality of  
$|G|+|A|+|B|$  implies that $B$  permutes  with $A\cap BR $. Also, $B$ permutes with 
$RA$ since $B(RA)=\langle A, B \rangle R$, so $AB=BA$ by Lemma 3.3, Part (i) and Theorem B.
 This contradiction shows that 
$A\leq BR$, so $BRA=BR$.  But  $R \leq O^{\sigma _{i}}(G) \leq N_{G}(B)$
 by Lemma 2.4(1), hence $B$ is normal in $BR$ and since $A \leq BR$ it follows that
$AB = BA$.   This contradiction shows that we 
have (3).
 
{\sl Final contradiction.}  Claims (2) and (3) imply that   $$G=\langle A, B 
\rangle O^{\sigma _{i}}(G)=\langle A, B 
\rangle\times  O^{\sigma _{i}}(G),$$ so every subgroup $H$ of $\langle A, B \rangle $ is
 $O^{\sigma _{i}}(G)$-invariant. It follows that every subgroup of   $\langle A, B 
\rangle $ is $\sigma$-permutable  in $G$ by Lemma 2.4(3) and Proposition 3.1. 
  Hence   ${\cal L}(\langle A, B 
\rangle)$   is a sublattice  of the distributive lattice   ${\cal L}.$   
 Thus $\langle A, B 
\rangle$ is cyclic by the Ore theorem  by \cite[Theorem 1.2.3]{R.}, so $AB=BA$, a
 contradiction.   The proposition is 
proved.

{\bf Corollary 3.5.} {\sl    If the   lattice ${\cal L}={\cal L}_{\sigma}(G)$  
 is distributive, then  every two members  $A$ and $B$ of ${\cal L}$  are  permutable.}
 
{\bf Proof.}   Suppose that this corollary is false and  let $G$ be  a counterexample
 with $|G|+|A|+|B|$  minimal.    

Let $R$ be a minimal normal subgroup of $G$. Lemma 3.2 implies that ${\cal 
L}_{\sigma}(G/R)$ is isomorphic to the interval $[G/R]$ in the modular 
lattice ${\cal L}_{\sigma per}(G)$. Therefore,  Lemma 3.2(2)   and
 the minimality of $G$ imply  that  
$(AR/R)(BR/R)=(BR/R)(AR/R)$. It follows that $RAB=\langle A, B \rangle R$
 is a subgroup of
 $G$,   so   $A_{G}=1$      and $B_{G}=1$. 
 Hence, because of Lemma 2.4(2), $A$ and $B$ are $\sigma$-nilpotent.  The 
minimality of $|G|+|A|+|B|$ implies that    for some  $i$ we 
have  $A, B \leq O_{\sigma _{i}}(G)$ and so $A, B \in {\cal L}_{\sigma _{i}}(G)$. 
  But  ${\cal L}_{\sigma _{i}}(G)$ 
is a sublattice of the distributive lattice ${\cal L}_{\sigma}(G)$ by Proposition 3.4(i).
 Therefore $AB=BA$ by Proposition 3.4(ii),  a contradiction. 
The corollary is proved.

{\bf Lemma 3.6} (See \cite[p. 59]{Gretz}). {\sl   A modular lattice $\cal L$ is  distributive
if and only if  
$ \cal L$  has no distinct elements $a, b$ and $c$ such that $ a \vee b=   a \vee c= 
  b \vee  c$ and $a  \wedge b= a  \wedge c=b  \wedge c$.  }

{\bf Lemma 3.7} (See \cite[Theorem 1.6.2]{R.}). {\sl Let $G=A\times B$,
 $f:A\to B$ is an isomorphism
 and $C=\{aa^{f}\ |\ \ a\in A\}$. Then $G=AC=BC$ and $A\cap C=1=B\cap C$.}

{\bf Lemma 3.8 } (See Corollary 2.4 and Lemma 2.5 in \cite{1}).  {\sl  The class  of all  
$\sigma$-nilpotent groups 
 ${\mathfrak{N}}_{\sigma}$      is closed under taking  
products of normal subgroups, homomorphic images and  subgroups.}

In view of Proposition 2.2.8  in \cite{15},  
 we get from Lemma 3.8 the following

{\bf Lemma 3.9.}    {\sl If 
$N$ is a normal subgroup of $G$, then
 $(G/N)^{{\frak{N}}_{\sigma}}=G^{{\frak{N}}_{\sigma}}N/N.$  }

{\bf Proof of Theorem A. }  {\sl Necessity.}  First note that every two members 
 of ${\cal L}$  are  permutable by Theorem B.  Moreover, since the lattice ${\cal 
L}_{n}(G)$  is a sublattice of the lattice $\cal L$,  it is distributive. 
 Since $G/D=G/G^{{\frak{N}}_{\sigma}}$ is $\sigma$-nilpotent by Lemmas 3.8 and 3.9,
 every subgroup $E$ of $G$ satisfying $D\leq E\leq G$ is
  $\sigma$-permutable 
in $G$ by Lemma 3.2(1) and Remark 1.2(ii). 
  Hence ${\cal L}(G/D)= {\cal L}_{\sigma per}(G/D)$
 is distributive and so 
$G/D$ is cyclic by the Ore theorem \cite[Theorem 1.2.3]{R.}.
 It is clear also that   $D$ 
is a  meet-distributive element of   ${\cal L}$.   Thus Conditions 
(i)-(iii) hold on $G$.

Now we show that Condition (iv) holds on $G$.  First note that  since,
 in view of Lemma 3.2,
 the lattice ${\cal L}_{\sigma per}(G/R)$ is  isomorphic
 to the interval $[G/R]$ in the distributive  
lattice ${\cal L}$, it is enough to consider the case when $\bar{G}=G$ and 
 $\bar{K}=K, \bar{H}=H, \bar{L}=L\in {\cal L}$. 

Suppose that $H\ne L$. Then $H\ne K$. Let $K < H_{0}\leq H$, where $H_{0}$ 
covers $K$ in $\cal L$, and let $L_{0}/K=(H_{0}/K)^{f}$, where 
$f: H/K   \to L/K$ is a $O^{\sigma _{i}}(G)$-isomorphism. For $g\in  
O^{\sigma _{i}}(G)$ and $l_{0}K=(hK)^{f}\in L_{0}/K$, where $h\in H$,  we have  
$$(l_{0}K)^{g}=((hK)^{f})^{g}=((hK)^{g})^{f}=(h^{g}K)^{f}=(h_{0}K)^{f},$$ 
where $h_{0}\in H_{0}$ since $H_{0}$ is  $O^{\sigma _{i}}(G)$-invariant by 
Lemma 2.4(1).  Hence  $(l_{0}K)^{g}\in L_{0}/K$. 
 It follows that   
 $L_{0}$ is  $O^{\sigma _{i}}(G)$-invariant 
 and so $L_{0}$ covers $K$ in 
$\cal L$ since the inverse map $f^{-1}:L/K   \to H/K$ is a $O^{\sigma _{i}}(G)$-isomorphism too.

 First assume that $H_{0}\ne L_{0}$  and let $E_{0}/K=\{hK(hK)^{f} | hK\in H_{0}/K\}$. 
Then  $$(H_{0}/K)(L_{0}/K)=(H_{0}/K)\times (L_{0}/K).$$ Indeed, if
 $H_0^x \not= H_0$ for some $x \in L_0$, then (i) and the fact that
 $H_0$ and $L_0$ cover $K$ in ${\cal L}$ would imply
that $\{K; H_0; H_0^x ; L_0; H_{0}L_{0}\}$ would be a diamond in the distributive
 lattice $\cal L$, contradicting Lemma 3.6.    
  Hence, by Lemma 3.7,  
 $E_{0}/K$ is a subgroup of $(H_{0}/K)\times (L_{0}/K)$ and  we have 
$$(H_{0}/K)\times (L_{0}/K)=(H_{0}/K)\times (E_{0}/K)=(L_{0}/K)\times 
(E_{0}/K).$$

Note that if $g\in  
O^{\sigma _{i}}(G)$ and $hK(hK)^{f}\in E_{0}/K$, then   
$$(hK(hK)^{f})^{g}=(hK)^{g}((hK)^{f})^{g}=(h^{g}K)(h^{g}K)^{f}\in E_{0}/K$$
 since $f_{H_{0}/K}$ is a 
$O^{\sigma _{i}}(G)$-isomorphism from $H_{0}/K$ onto $L_{0}/K=(H_{0}/K)^{f}$.
Hence 
$E_{0}/K$ is $O^{\sigma _{i}}(G)$-invariant, so $O^{\sigma _{i}}(G)\leq N_{G}(E_{0})$.
  Therefore $H_{0}$, $L_{0}$  
and   $E_{0}$ are distinct elements of $\cal L$  such that  $H_{0}\cap L_{0}= 
H_{0}\cap E_{0}=L_{0}\cap E_{0}=K $ and   $H_{0}L_{0}= 
H_{0}E_{0}=L_{0}E_{0} $, which is impossible by Lemma 3.6 since 
$H_{0}L_{0}$ is a
 $\sigma$-permutable 
subgroup of $G$.   
Therefore  $H_{0}= L_{0}$. Now $f$ induces a $O^{\sigma _{i}}(G)$-isomorphism $f':H/H_{0}\to 
L/H_{0}$ and an obvious induction yields that $H=L$. 
Hence we  have (iv).

{\sl Sufficiency.}  This follows from the following

{\bf Proposition 3.10.} {\sl   Let $D=G^{{\frak{N}}_{\sigma}} $ and 
 ${\cal L}={\cal L}_{\sigma per}(G)$. 
Suppose that  
 the following conditions hold:}

(i) {\sl    Every two members 
 of ${\cal L}$  are  permutable.  }

(ii) {\sl The lattice ${\cal L}_{n}(G)$ of all  normal subgroups of $G$ is distributive. }

(iii) {\sl $G/D$ is cyclic and $D$ is a  meet-distributive element of   ${\cal L}$.}

(iv) {\sl  In every factor group $\bar{G}=G/R$,   any 
 two $O^{\sigma _{i}}(\bar{G})$-isomorphic sections
 $\bar{H}/\bar{K}$ and $\bar{L}/\bar{K}$, where   $\bar{K}, \bar{H}, 
\bar{L}\in {\cal L}_{\sigma _{i}}(\bar{G})$ (for some $i$) and the subgroups
 $\bar{H}$ and $\bar{L}$
 cover $\bar{K}$
 in ${\cal L}_{\sigma  per}(\bar{G})$,  coincide. }

{\sl  Then $\cal L$ is distributive. }

{\bf Proof.}  
Suppose that this is false and let $G$ be a counterexample of minimal 
order.    

First note that if  $A, B, C\in  {\cal L}_{\sigma per}(G)$ and $A\leq C$, then  $$C\cap 
\langle A, B\rangle=C\cap AB=A(C\cap B)=\langle A, C\cap B\rangle$$ by Condition (i),  so 
 the lattice ${\cal L}_{\sigma}(G)$ is modular.   Hence, by Lemma 3.6,
 there are 
  distinct $\sigma$-permutable 
 subgroups $A$, $B$ and
 $C$ of $G$ such 
that  for some  $\sigma$-permutable  subgroups $E$ and $T$ of $G$ we
 have $E=A\cap B=A\cap C=B\cap C$
   and $T=AB=AC=BC$.  

(1) {\sl The lattice ${\cal L}_{\sigma per}(G/R)$ is distributive for each non-identity
 normal subgroup $R$ of $G$.  }

In view of the choice of $G$, it is enough to show that Conditions  
(i), (ii),  (iii) and (iv) hold for $G/R$. 

  Let $K/R, H/R  \in {\cal L}_{\sigma per}(G/R)$. Then $K, H \in {\cal 
L}(G)$ by Lemma 3.2(1)  and  so $KH=HK$ by Condition (i),
 which implies that $(K/R)(H/R)=(H/R)(K/R)$.   
  It is clear also that  the lattice ${\cal L}_{n}(G/R)$ is isomorphic to
 some sublattice of the lattice ${\cal L}_{n}(G)$, so ${\cal L}_{n}(G/R)$  is distributive.
Thus Conditions  
(i)  and  (ii)  hold on $G/R$.

By Lemma 3.9 we have
 $$(G/R)^{{\frak{N}}_{\sigma}}=G^{{\frak{N}}_{\sigma}}R/R=DR/R.$$ Thus 
$$(G/R)/(G/R)^{{\frak{N}}_{\sigma}}=(G/R)/(DR/R)  \simeq G/DR\simeq   
(G/D)/(GR/D)$$  is cyclic by  Condition (iii).
Conditions (i) and (iii) imply  that $$D\cap \langle K, H\rangle =D\cap 
\langle D\cap K, D\cap H\rangle= 
   (D\cap K)(D\cap H)$$ since  $D\cap K$ and $D\cap H$
 are $\sigma$-permutable in $G$ by Theorem B,  so  
   $$(G/R)^{{\frak{N}}_{\sigma}}\cap \langle (K/R), (H/R)\rangle = (DR/R)\cap 
(K/R)(H/R)=(DR\cap KH)/R=R(D\cap KH)/R=$$$$=R(D\cap K)(D\cap H)/R=
((D\cap K)R/R)((D\cap H)R/R)=$$$$=((DR/R)\cap (K/R))((DR/R)\cap (H/R))=
\langle(G/R)^{{\frak{N}}_{\sigma}}\cap (K/R)),   (G/R)^{{\frak{N}}_{\sigma}}\cap
 (H/R)\rangle.$$ 
  Hence $(G/R)^{{\frak{N}}_{\sigma}}$ is a  meet-distributive element of   ${\cal 
L}_{\sigma per}(G/R)$.  Thus    Condition  
  (iii)  hold on $G/R$.    
Finally,  Condition
 (iv), evidently,   hold on $G/R$. Thus we have (1).

 (2) $E_{G}=1$ (In view of Lemma 3.2, this follows from Claim (1), 
Lemma 3.6 
and the choice of $G$).

(3)  {\sl $A_{G}B_{G}\cap A_{G}C_{G}\cap B_{G}C_{G}=1$}.

Since $A\cap B=E$,  we have $B_{G}\cap A_{G}\leq E_{G}=1$ by 
  Claim (2). Similarly, $B_{G}\cap C_{G}=1$   and $A_{G}\cap C_{G}=1$. Therefore   

 $$(A_{G}B_{G}\cap A_{G}C_{G})\cap 
B_{G}C_{G}=
A_{G}(B_{G}\cap A_{G}C_{G})\cap B_{G}C_{G}=$$$$=
A_{G}(B_{G}\cap A_{G})(B_{G}\cap C_{G})\cap B_{G}C_{G}=$$ $$=A_{G}\cap B_{G}C_{G}=
(A_{G}\cap B_{G})(A_{G}\cap C_{G})=1$$   by Claim (ii).

(4) {\sl The subgroup $T$ is $\sigma$-nilpotent. }

Note that 
$$T/A_{G}B_{G}=AB/A_{G}B_{G}
=(AA_{G}B_{G}/A_{G}B_{G})(BA_{G}B_{G}/A_{G}B_{G}),$$ where 
$$AA_{G}B_{G}/A_{G}B_{G}\simeq A/A\cap
 A_{G}B_{G}=
$$$$=A/A_{G}(A\cap B_{G})\simeq 
(A/A_{G})/(A_{G}(A\cap B_{G})/A_{G})$$  and 
 $$BA_{G}B_{G}/A_{G}B_{G}\simeq  
(B/B_{G})/(B_{G}(B\cap A_{G})/B_{G})$$    are $\sigma$-nilpotent by   
Lemma 2.4(2). The subgroups  $AA_{G}B_{G}/A_{G}B_{G}$ and $BA_{G}B_{G}/A_{G}B_{G}$  are
 $\sigma$-subnormal in $G/A_{G}B_{G}$ by Lemmas   2.2 and 3.2(2).
  Hence $T/A_{G}B_{G}$ is $\sigma$-nilpotent by Lemma 2.3(6).
  Similarly,   $T/A_{G}C_{G}$ and 
$T/C_{G}B_{G}$  are $\sigma$-nilpotent.  Hence from Claim  (3) it 
follows  that $T\simeq T/(A_{G}B_{G}\cap A_{G}C_{G}\cap B_{G}C_{G})$  is 
$\sigma$-nilpotent by Lemma 3.8.

 (5)  {\sl For some prime $i$, there are distinct
 $\sigma _{i}$-subgroups $A_{i}, B_{i}, C_{i}\in {\cal L}$ such that 
 $H_{i}=A_{i}B_{i}=A_{i}C_{i}=B_{i}C_{i}$ and 
 $K_{i}=A_{i}\cap B_{i}=A_{i}\cap C_{i}=B_{i}\cap C_{i}$ are $\sigma$-permutable
 subgroups of $G$. }

Let   $\sigma _{i}\in \sigma (T)$, that is,  $\sigma _{i}\cap \pi (T)\ne \emptyset$. Then, by Claim (4), $H_{i}=O_{\sigma _{i}}(T)$ is the  Hall
 $\sigma _{i}$-subgroup of $T$ and  $A_{i}=O_{\sigma _{i}}(A)$,
 $B_{i}=O_{\sigma _{i}}(B)$ and $C_{i}=O_{\sigma _{i}}(C)$  are  the  Hall
 $\sigma _{i}$-subgroups of 
$A$, $B$ and $C$, respectively. Hence 
$H_{i}=A_{i}B_{i}=A_{i}C_{i}=B_{i}C_{i}$. Moreover, $A_{i}$, $B_{i}$ and 
$C_{i}$ are $\sigma$-permutable  in $G$ by Theorem C.
It is clear also  that 
  $K_{i}=A_{i}\cap B_{i}=A_{i}\cap C_{i}=B_{i}\cap C_{i}=O_{\sigma _{i}}(E)$.

Suppose that $A_{i}= B_{i}$.
 Then $H_{i}=A_{i}B_{i}=A_{i}=B_{i}= K_{i}\leq C_{i}\leq H_{i}$. 
Hence $A_{i}= B_{i}= C_{i}$. Therefore,  since $A\ne B\ne C$ and $A\ne  C$, there is
 $\sigma _{i}\in \sigma (T)$ such that $A_{i}\ne B_{i}\ne C_{i}$ and $A_{i}\ne  C_{i}$. 
Finally,  $H_{i}$ and $K_{i}$   are $\sigma$-permutable subgroups of $G$ by Condition (i),  
 so we have (5).

(6)   {\sl There are distinct  $\sigma _{i}$-subgroups $A_{0}, B_{0}, C_{0}\in {\cal L}$
 such that 
 $H_{0}=A_{0}B_{0}=A_{0}C_{0}=B_{0}C_{0}$ and 
 $K_{0}=A_{0}\cap B_{0}=A_{0}\cap C_{0}=B_{0}\cap C_{0}$ are
 $\sigma$-permutable subgroups of $G$ 
 and   
$A_{0}, B_{0}, C_{0}$ are normal subgroups of $O^{\sigma _{i}}(G)$.  
}

Let $A_{0}=A_{i}\cap D$,  $B_{0}=B_{i}\cap D$  and $C_{0}=C_{i}\cap D$.
Then    $A_{0}$, $B_{0}$ and 
$C_{0}$ are $\sigma$-permutable $\sigma _{i}$-subgroups of $G$ by Claim  (5) and Theorem B. 
Moreover, Claim (5) implies that 
 $$K_{0}=A_{0}\cap B_{0}= A_{i}\cap B_{i}\cap D=A_{i}\cap C_{i}\cap D=A_{0}\cap C_{0} =
 B_{i}\cap C_{i}\cap D=B_{0}\cap C_{0}.$$   
Since $D$ is a    meet-distributive element of   ${\cal L}$ by Condition (iii),
 $$H_{0}=D\cap A_{i}B_{i}=(D\cap A_{i})(D\cap B_{i})=A_{0}B_{0}=A_{0}C_{0}=D\cap A_{i}C_{i}=
D\cap B_{i}C_{i}=B_{0}C_{0}.$$

Now   we show that $A_{0}, B_{0}, C_{0}$ are distinct elements of 
${\cal L}$. First note that  $$|H_{i}:K_{i}|=|A_{i}:K_{i}||B_{i}:K_{i}|=
|A_{i}:K_{i}||C_{i}:K_{i}|=|B_{i}:K_{i}||C_{i}:K_{i}|,$$ 
so $|A_{i}:K_{i}|=|B_{i}:K_{i}|=|C_{i}:K_{i}|$. Hence $|A_{i}|=|B_{i}|=|C_{i}|.$
Suppose that  $A_{0}=B_{0}$.
Then
 $$D\cap H_{i}=D\cap A_{i}B_{i}=(D\cap A_{i})(D\cap B_{i})=
A_{0}B_{0}=A_{0}=B_{0}=D\cap K_{i}.$$ 
Hence $K_{i}D\cap H_{i}=K_{i}(D\cap H_{i})=K_{i}$ is normal in $H_{i}$ and
 $$DH_{i}/DK_{i}\simeq H_{i}/(H_{i}\cap 
K_{i}D)=H_{i}/K_{i}(H_{i}\cap D)=H_{i}/K_{i}$$ is cyclic since $G/D$ is cyclic by
 Condition (iii). 
On the other hand, $H_{i}/K_{i}=(A_{i}/K_{i})(B_{i}/K_{i})$, where
  $|A_{i}/K_{i}|=|B_{i}/K_{i}|$, so $A_{i}/K_{i}=B_{i}/K_{i}=1$,
 which implies that $A_{i}=B_{i}$.  
This contradiction shows that  $A_{0}\ne B_{0}$. Similarly, 
one can show that   $A_{0}\ne C_{0}$  and $B_{0}\ne C_{0}$. 
 Finally, $A_{0}, B_{0}, C_{0}$ are 
normal subgroups of $O^{\sigma _{i}}(G)$ by Lemma 2.4(1), and Claim (5) 
and Theorem B imply 
 that $K_{0}$ and $H_{0}$ are $\sigma$-permutable in $G$.

(7) {\sl $A_{0}/K_{0}$ and $B_{0}/K_{0}$ are $O^{\sigma _{i}}(G)$-isomorphic.}

From  Claim (6) we get that  $$H_{0}/K_{0}=(A_{0}/K_{0})\times (B_{0}/K_{0})= 
(A_{0}/K_{0})\times (C_{0}/K_{0})=(B_{0}/K_{0})\times (C_{0}/K_{0}).$$  Therefore   
$$A_{0}/K_{0}\simeq ((A_{0}/K_{0})\times
 (C_{0}/K_{0}))/(C_{0}/K_{0})=(H_{0}/K_{0})/(C_{0}/K_{0})$$ and
 $$B_{0}/K_{0}\simeq ((B_{0}/K_{0})\times
 (C_{0}/K_{0}))/(C_{0}/K_{0})=(H_{0}/K_{0})/(C_{0}/K_{0})$$ 
  are $O^{\sigma _{i}}(G)$-isomorphisms 
by Lemma 2.4(1).  Hence we have (7).

{\sl Final contradiction.}   Let $f: A_{0}/K_{0}  \to B_{0}/K_{0}$ be a   
$O^{\sigma _{i}}(G)$-isomorphism. Let $K_{0} < X\leq A_{0}$, where  $X$  
covers  $K_{0}$ in ${\cal L}$.  Then  $X/K_{0}$ is a chief factor of  
$O^{\sigma _{i}}(G)$ by Lemma 2.4(3)  and Proposition 3.1, so $L/K_{0}=f(A_{0}/K_{0})$  is also a chief factor of  
$O^{p}(G)$.  Hence $L$  
covers  $K_{0}$ in ${\cal L}$  by Lemma 2.4(1). Now $f$ induces a 
$O^{\sigma _{i}}(G)$-isomorphism from $X/K_{0}$ onto $L/K_{0}$ and  so
 $L=T$ by Condition (iv). Hence $K_{0} < A_{0}\cap B_{0}$, contrary to (6).

 The proposition is proved.

\end{document}